\documentstyle[12pt]{article}
\begin{document}
\newtheorem{Def}{Definition}
\newtheorem{thm}{Theorem}
\newtheorem{lem}{Lemma}
\newtheorem{rem}{Remark}
\newtheorem{prop}{Proposition}
\newtheorem{cor}{Corollary}
\newtheorem{clm}{Claim}
\title
{A fully nonlinear version of the Yamabe problem and a Harnack type
inequality}
\author{Aobing Li \ \ \ \ \& \ \ \ YanYan Li
\\
Department of Mathematics\\
Rutgers University\\
110 Frelinghuysen Rd.\\
Piscataway, NJ 08854
\\
emails:  \  aobingli@math.rutgers.edu, \
yyli@math.rutgers.edu
}
\input { amssym.def}
\date{}
\maketitle

We present some results in
 \cite{LL2}, a continuation
of our earlier works \cite{LL0} and \cite{LL}.
One result is the  existence and compactness
of solutions to a fully nonlinear version of the Yamabe problem
on locally conformally flat Riemannian manifolds, and the other
is a Harnack type inequality for general
conformally invariant fully nonlinear second order elliptic
equations.

Let $(M,g)$ be an $n-$dimensional, compact, smooth
Riemannian manifold without boundary, $n\ge3$, 
 consider the
Weyl-Schouten
 tensor
$
A_g=\frac{1}{n-2}\left(Ric_g-\frac{R_g}{2(n-1)}g\right),
$
where $Ric_g$ and $R_g$ denote respectively the
Ricci tensor and the scalar curvature associated with $g$.
We use $\lambda(A_g)$ to denote
the eigenvalues of $A_g$ with respect to $g$.

Let $\hat g=u^{\frac 4 {n-2} }g$ be a conformal change of
metrics, then (see, e.g., \cite{V3}),
\begin{equation}
A_{\hat g}=-\frac{2}{n-2}  u^{-1}\nabla_{g}^2u+ \frac{2n}{(n-2)^2} u^{-2}
\nabla_{g}u \otimes\nabla_{g}u -\frac{2}{(n-2)^2} u^{-2}
|\nabla_{g}u|_{g}^2g+A_{g}.
\label{change}
\end{equation}

Let
\begin{equation}
\Gamma\subset \Bbb R^n \ \mbox{be
an open convex cone  with  vertex at the origin },
\label{gamma0}
\end{equation}
\begin{equation}
\{\lambda\in \Bbb R^n\ |\ \lambda_i>0, 1\le i\le n\}
\subset \Gamma\subset \{\lambda\in \Bbb R^n\ |\ \sum_{i=1}^n  \lambda_i>0\},
\label{gamma1}
\end{equation}
\begin{equation}
\Gamma\ \mbox{is symmetric in the}\ \lambda_i,
\label{gamma1new}
\end{equation}
\begin{equation}
f\in C^\infty(\Gamma)\cap
C^0(\overline \Gamma)\
\mbox{be  concave and symmetric in the}\ \lambda_i,
\label{hypo0}
\end{equation}
\begin{equation}
f=0\ \mbox{on}\ \partial \Gamma;\quad
f_{\lambda_i}>0\ \mbox{on}\ \Gamma\ \forall\  1\le i\le n,
\label{aa5}
\end{equation}
\begin{equation}
\lim_{s\to\infty}f(s\lambda)=\infty,\qquad
\forall\ \lambda\in \Gamma.
\label{aa6}
\end{equation}

\begin{thm} \ (\cite{LL2})\ For  $n\ge 3$, let  $(f,\Gamma)$ satisfy
(\ref{gamma0}), (\ref{gamma1}), (\ref{gamma1new}),
 (\ref{hypo0}), (\ref{aa5})
and (\ref{aa6}), and let $(M, g)$ be an
$n-$dimensional  smooth  compact locally conformally flat
 Riemannian  manifold without boundary satisfying
\begin{equation}
\lambda\left(A_{g}\right)\in \Gamma,\qquad \mbox{on}\ M.
\label{aa1}
\end{equation}
Then there exists some smooth positive function $u$ on $M$ such that
$\hat g=u^{\frac 4{n-2}} g$ satisfies
\begin{equation}
f\left(\lambda(A_{\hat g})\right)=1,
\quad \lambda(A_{\hat g})\in \Gamma, \qquad \mbox{on}\ M.
\label{ghat2}
\end{equation}
Moreover, if  $(M,g)$ is not conformally diffeomorphic
to the  standard $n-$sphere,  all solutions of the above
satisfy, for any positive integer $m$, that
\begin{equation}
\|u\|_{C^m(M,g)}+\|u^{-1}\|_{C^m(M,g)}\le C,
\label{aa7}
\end{equation}
where $C$ is some  constant
depending  only on $(M,g)$,  $(f,\Gamma)$ and
$m$.
\label{t1}
\end{thm}
\begin{rem} In the proof of Theorem \ref{t1},, we see
that the $C^0$ and $C^1$ apriori estimates
above do not require the concavity of $f$.
More precisely,  without the concavity assumption on
$f$ in the statement of Theorem \ref{t1}, and
when  $(M,g)$ is not conformally diffeomorphic
to the  standard $n-$sphere,  all  solutions
of (\ref{ghat2}) satisfy, for some constant
$C$ depending only on $(M,g), b, \delta_1$ and $\delta_2$, that
$
\|u\|_{ C^1(M,g) }+
\|u^{-1}\|_{ C^1(M,g) }\le C.
$
\label{rem1-1}
\end{rem}

For  $1\le k\le n$,
let
$
\sigma_k(\lambda)=\sum_{1\le i_1<\cdots<i_k\le n}\lambda_{i_1}
\cdots \lambda_{i_k}, \lambda=(\lambda_1,
 \cdots, \lambda_n)\in \Bbb R^n,
$
denote the $k-$th symmetric function, and let
$\Gamma_k$ denote the
connected component of $\{\lambda\in \Bbb R^n\ |\
\sigma_k(\lambda)>0\}$ containing  the
positive cone $\{\lambda\in \Bbb R^n\ |\
\lambda_1, \cdots, \lambda_n>0\}$.
Then (see \cite{CNS})
$(f, \Gamma)=(\sigma_k^{\frac 1k}, \Gamma_k)$ satisfies
the hypothesis in Theorem \ref{t1}.

For $(f,\Gamma)=(\sigma_1, \Gamma_1)$,
hypothesis (\ref{aa1})  
is equivalent to $R_g>0$ on $M$, and therefore 
Theorem \ref{t1} in this case
 is the Yamabe
problem for  locally conformally flat manifolds
with positive Yamabe invariants,
 and  the result  is due to Schoen (\cite{S0} and \cite{S1}).
The Yamabe conjecture was proved through the work of Yamabe,
Trudinger, 
Aubin and Schoen.
For  $(f,\Gamma)=(\sigma_k^{\frac 1k},
\Gamma_k)$ with $k=2$ and $n=4$,  the result was proved without the
locally conformally flatness hypothesis of the manifold by
Chang, Gursky and Yang \cite{CGY2}.
For  $(f,\Gamma)=(\sigma_k^{\frac 1k},
\Gamma_k)$ with $k=n\ge 3$,  some existence
result was established by Viaclovsky \cite{V1} for a class of manifolds
which are not necessarily locally conformally flat.
For  $(f,\Gamma)=(\sigma_k^{\frac 1k},
\Gamma_k)$, $n\ge 3$, $1\le k\le n$, Theorem \ref{t1}
was established in \cite{LL0} and \cite{LL}; while 
the
existence part in the case  $k\neq \frac n2$
 was independently  obtained by
Guan and Wang in \cite{GW2}.
Subsequently,  Guan,  Viaclovsky and
Wang \cite{GVW}  proved the algebraic fact 
that  $\lambda(A_g)\in \Gamma_k$ for
$k\ge \frac n2$ implies  the positivity of the  Ricci tensor,
and therefore both the existence
 and compactness results in this case follow from known results.
More recently, Gursky and Viaclovsky \cite{GV2}
have obtained    
 existence results
for $(f,\Gamma)=(\sigma_k^{\frac 1k}, \Gamma_k)$, $n=3,4$,
on general Riemannian manifolds.

A Liouville type theorem 
for $(f,\Gamma)=(\sigma_k^{\frac 1k}, \Gamma_k)$
was established in \cite{LL}.
  The crucial ingredient in our proof of the
Liouville type theorem is a Harnack type inequality 
for $(f,\Gamma)=(\sigma_k^{\frac 1k}, \Gamma_k)$ established in the
same paper.   In \cite{LL2}, we have established
 the Harnack type inequality
for general
conformally invariant fully nonlinear second order elliptic
equations.
In the following, ${\cal S}^{n\times n}$
denotes the set of $n\times n$ real symmetric matrices,
${\cal S}^{n\times n}_+\subset {\cal S}^{n\times n}$ denotes the set
of positive definite matrices,
 $O(n)$ denotes the set of $n\times n$ real orthogonal matrices,
and $I$ denotes the $n\times n$ identity matrix..
It was show in \cite{LL} that
$H(\cdot, u, \nabla u, \nabla^2 u)$
is  conformally invariant on $\Bbb R^n$ 
(see \cite{LL} for the definition)
if and only if
$
 H(\cdot, u, \nabla u, \nabla^2 u)
\equiv F(A^u),
$
where
\begin{equation}
A^u:= -\frac{2}{n-2}u^{  -\frac {n+2}{n-2} }
\nabla^2u+ \frac{2n}{(n-2)^2}u^ { -\frac {2n}{n-2} }
\nabla u\otimes\nabla u-\frac{2}{(n-2)^2} u^ { -\frac {2n}{n-2} }
|\nabla u|^2I,
\label{2}
\end{equation}
and $F$  is invariant under orthogonal conjugation,
i.e.,
\begin{equation}
F(O^{-1}MO)=F(M), \qquad \forall\ M\in {\cal S}^{n\times n},\
\forall\ O\in O(n).
\label{F1}
\end{equation}

Let
 $U\subset {\cal S}^{n\times n}$ 
be  an open set satisfying
\begin{equation}
O^{-1}UO=U,\qquad\forall\ O\in O(n),
\label{U1}
\end{equation}
\begin{equation}
U\cap \{M+tN\ |\ 0<t<\infty\}\ \mbox{is convex}\qquad \forall\
M\in {\cal S}^{n\times n}, N\in {\cal S}^{n\times n}_+.
\label{U2}
\end{equation}

Let $F\in C^1(U)$ satisfy (\ref{F1})
and
\begin{equation}
\left(F_{ij}(M)\right)>0,\qquad \forall\ M\in U,
\label{F2}
\end{equation}
where $F_{ij}(M):=\frac{\partial F}{ \partial M_{ij} }(M)$.
We assume that for some $\delta>0$,
\begin{equation}
F(M)\neq 1\qquad\forall\
M\in U\cap \{M\in  {\cal S}^{n\times n}\
|\ \|M\|:=(\sum_{i,j}M_{ij}^2)^{\frac 12} <\delta\}.
\label{F3}
\end{equation}

For
$
F_k(M):=\sigma_k^{ \frac 1k}(\lambda(M)),
$ and $U_k:=\{M\in {\cal S}^{n\times n}\
|\ \lambda(M)\in \Gamma_k\},
$
it is well known that  $(F,U)=(F_k, U_k)$ satisfies  
(\ref{U1}), (\ref{U2}), (\ref{F1}) and (\ref{F3}).
In the following we use $B_R$  to denote a
ball in $\Bbb R^n$ which is
of radius $R$ and centered at the origin.
\begin{thm}\ (\cite{LL2})\
For $n\ge 3$, let $U\subset {\cal S}^{n\times n}$ satisfy
(\ref{U1}) and (\ref{U2}), and let $F\in C^1(U)$ satisfy
(\ref{F1}), (\ref{F2}) and  (\ref{F3}).
For $R>0$, let $u\in C^2(B_{3R})$ be a positive solution of 
\begin{equation}
F(A^u)=1,\quad A^u\in U,\quad \mbox{in}\quad B_{3R}.
\label{t2e1}
\end{equation}
Then
\begin{equation}
(\sup_{B_R}u)(\inf_{B_{2R}}u)\le C(n)\delta^{ \frac {2-n}2 }R^{2-n},
\label{t2e3}
\end{equation}
where $C(n)$ is some constant depending  only on $n$.
\label{t2}
\end{thm}

\begin{rem}  In Theorem \ref{t2},
 there is no concavity assumption 
on $F$ and the constant $C(n)$ can be given explicitly.
\end{rem}

\begin{rem}  The Harnack type inequality 
(\ref{t2e3}) for
$(F,U)=(F_1, U_1)$ was obtained by Schoen in \cite{S2} based
on a Liouville type theorem of Caffarelli, Gidas and Spruck
in \cite{CGS}. 
Li and
Zhang gave in \cite{LZ} a different proof
of Schoen$'$s   Harnack type inequality without using the
Liouville type theorem.
For $(F,U)=(F_k, U_k)$, $1\le k\le n$, the Harnack type inequality
was established in our earlier work  \cite{LL}.
There are two new ingredients in our proof of
Theorem \ref{t2}.  One is that we have developed, along the line of
\cite{LL}, new $C^0$ and $C^1$ estimates which
allow us to extend the Harnack type
inequality in \cite{LL} to this generality, and the other is that 
we have given a direct proof which makes it possible to
give an explicit constant $C$ in (\ref{t2e3}).
Arguments  in  \cite{S2}, \cite{LZ} and \cite{LL}
were indirect and therefore
no explicit value of $C$ was available, even in the case 
$(F,U)=(F_1, U_1)$.
\end{rem}

We first present our proof of
Theorem \ref{t1}, more 
details can be
found in \cite{LL2}.
As explained in \cite{LL2}, we may further assume without loss of generality
that $f$ is homogeneous of degree $1$.
By (\ref{aa5}) and (\ref{aa6}), there exists a unique $b>0$ such that
$f(be)=1,
$
where $e=(1,\cdots, 1)$.
By (\ref{aa5}), there exists some $\delta_1>0$ such that
\begin{equation}
f(\lambda)<1,\qquad \forall\ \lambda\in \Gamma, |\lambda|<\delta_1.
\label{9new}
\end{equation}

Fix some constant $\delta_2$ such that
\begin{equation}
0<\delta_2\le \min_{  x\in M} f(\lambda(A_g(x)).
\label{delta2}
\end{equation}

Let $(\widetilde M,
\widetilde g)$ denote the universal cover of $(M,g)$,
 with $i:\widetilde M\to M$
a covering map and $ \widetilde g=i^*g$.
By a  theorem of Schoen and Yau in \cite{SY1},
there exists an injective  conformal immersion
$
\Phi: (\widetilde M,  \widetilde g) \to (\Bbb S^n, g_0),
$
where $ g_0$ denotes the standard metric on $\Bbb S^n$.
Moreover, $
\Omega:=\Phi(\widetilde M)
$ is either $\Bbb S^n$ or   an open and dense
subset of $\Bbb S^n$.
Fix a compact subset $E$ of $\widetilde M$ such that
$i(E)=M$.

To prove Theorem \ref{t1}, we will establish (\ref{aa7}) first.
 In the following, let
  $u\in C^\infty(M)$ be a  positive
solution of (\ref{ghat2}) with $\hat g=u^{  \frac 4{n-2} }g$.

\noindent {\bf Step 1.}\  For some  positive constant  $C$
 depending only on $(M,g)$, $b$,  $\delta_1$ and $\delta_2$, we have
\begin{equation}
\frac 1C \le u\le C,\quad |\nabla_g u|\le C\qquad\mbox{on}\ M.
\label{estimate}
\end{equation}

We will use notation
$
F(A_g):= f(\lambda(A_g)).
$
We distinguish into two cases.

\noindent {\bf Case 1.}\
$
\Omega= \Bbb S^n;
\qquad$
 {\bf Case 2.}\
$
\Omega\neq \Bbb S^n.
$

 In Case 1,
$(\Phi^{-1})^*\widetilde g=\eta^{\frac 4{n-2}}g_0$ on
$\Bbb S^n$, where $\eta$ is a positive smooth function on
$\Bbb S^n$.
Let $\tilde u= u\circ i$.  Since
$F\left(A_{ \tilde u^{\frac 4{n-2}}
\widetilde g }
\right) =1$ on $\widetilde M$, we have
$
F\left(
A_{ [(\tilde u\circ \Phi^{-1})\eta]^{\frac 4{n-2}}g_0 }\right)=1,
$ on $\Bbb S^n.
$
By corollary 1.1 in \cite{LL},
$ (\tilde u\circ \Phi^{-1})\eta=a |J_{ \varphi}|^{ \frac {n-2}{2n} }$
for some positive constant $a$ and some conformal diffeomorphism
$\varphi: \Bbb S^n\to \Bbb S^n$.  Since
$\varphi^* g_0=|J_{ \varphi}|^{ \frac 2n}g_0$,
we have, by the above equation, that
$
f(a^{-\frac 4{n-2}}(n-1)e)=
f(a^{-\frac 4{n-2}}\lambda(A_{g_0}))=1,
$
i.e. $(n-1)a^{-\frac 4{n-2}}=b$.
With this explicit formula of $\tilde u$, estimate
(\ref{aa7}) can be established without much difficulty.

In Case 2,    $(\Phi^{-1})^*\widetilde g=\eta^{\frac 4{n-2}}g_0$
on $\Omega$ where,  by 
\cite{SY1},  $\eta$ is a positive
smooth function in $\Omega$ satisfying
$\lim_{z\to \partial\Omega}\eta(z)=\infty$.
Recall that $
\Omega
$ is   an open and dense
subset of $\Bbb S^n$.
Let $u(x)=\max_M u$ for some $x\in M$, and let
$i(\tilde x)=x$ for some $\tilde x\in E$.
By composing with a rotation of $\Bbb S^n$, we
may assume without loss of generality that
$\Phi(\tilde x)=S$, the south pole of $\Bbb S^n$.
Let $P: \Bbb S^n \to \Bbb R^n$ be the stereographic projection, and
let $v$ be  the positive function on the open subset
 $P(\Omega)$ of $\Bbb R^n$
determined by
$(P^{-1})^*(\eta^{ \frac 4{n-2} }g_0)=v^{ \frac 4{n-2} }g_{flat}$,
where $g_{flat}$ denotes the
Euclidean metric on $\Bbb R^n$.
Then for some $\epsilon>0$, depending only on $(M,g)$ ,
we have
$
B_{9\epsilon}:=\{x\in \Bbb R^n\
|\ |x|<9\epsilon\}\subset P(\Omega),
$
and
$
dist_{ g_{flat}} \big(P(\Phi(E)),
\partial  P(\Omega)\big)>9\epsilon.
$
Let
$
\hat u=(\tilde u\circ \Phi^{-1}\circ P^{-1})v$ on
$
P(\Omega),
$
we have, by (\ref{change}),
$
f\big(\lambda(A^{\hat u})\big)=1$ and
$\lambda(A^{\hat u})\in \Gamma.
$
By the property of $\eta$, we know that
\begin{equation}
\lim_{y\to \bar y,
y\in P(\Omega)}\hat u(y)= \infty\quad
\forall \ \bar y\in \partial
P(\Omega),
\label{unbounded}
\end{equation}
and, if the north pole of $\Bbb S^n$  does not
belong to $\Omega$,
\begin{equation}
\lim_{ y\in P(\Omega), |y|\to\infty}
(|y|^{n-2} \hat u(y))=\infty.
\label{20new}
\end{equation}
 
By  a moving sphere argument
(i.e. moving plane method together with conformal invariance of the
equation) as in \cite{LL2}, 
we have,
 for  every $x\in \Bbb R^n$ satisfying
$dist_{g_{flat}}(x, P(\Phi(E)))<2\epsilon$, that
\begin{eqnarray}
\hat u_{x,\lambda}(y):=&&\frac { \lambda^{n-2} }
{ |y-x|^{n-2} } \hat u( \frac { \lambda^2 (y-x) }{ |y-x|^2 } )\le
\hat u(y),\nonumber\\
&&\qquad\qquad\qquad
\forall\ 0<\lambda<4\epsilon,
|y-x|\ge\lambda,~y\in P(\Omega).
\label{move}
\end{eqnarray}

The following calculus lemma is established
in \cite{LL2}.
\begin{lem} Let $a>0$ be a constant and
let  $B_{8a}\subset \Bbb R^n$ be the ball of radius $8a$ and centered at the
origin, $n\ge 3$.
Assume that $u\in C^1(B_{8a})$ is a non-negative function
satisfying
$$
u_{x,\lambda}(y)\le u(y),\quad\forall \
x\in B_{4a}, \ y\in B_{8a},~0<\lambda<2a,~\lambda<|y-x|,
$$
where
$u_{x,\lambda}(y):=\Big(\frac{\lambda}{|y|}\Big)^{n-2}u\Big(x+\frac
{\lambda^2(y-x)}{|y-x|^2}\Big)$.
Then
\[
|\nabla u(x)|\le  \frac {n-2}{2a} u(x),\quad\forall |x|<a.
\]
\label{l2}
\end{lem}

By (\ref{move}) and the above lemma, we have
$
|\nabla(\log \hat u)(y)|\le C(\epsilon)$,
\newline
$
\forall\ dist_{g_{flat}}(y, P(\Phi(E)))< \epsilon.
$
Thus, for some positive constant  $C$ depending only on $(M,g)$,
$
|\nabla_g \log u|\le C$ on $ M,
$
and
\begin{equation}
\sup_{ B_{\epsilon} }\hat u\le C \inf_{ B_{\epsilon} }\hat u.
\label{harnack}
\end{equation}

Let  $\beta>0$ be the constant such  that
$
\xi(y):=\beta(\epsilon^2-|y|^2)
$
has the property that
$
\hat u\ge \xi$ on $ B_\epsilon,
$
and, for some $\bar y\in B_\epsilon$,
$
\hat u(\bar y)=\xi(\bar y).
$
It follows that
$
\nabla \hat u(\bar y)= \nabla \xi(\bar y),
\quad (D^2\hat u(\bar y))\ge  (D^2 \xi(\bar y)),
$ and $ 
A^{ \hat u}(\bar y)\le A^{\xi}(\bar y).
$
By  (\ref{harnack}) and the definition of $\xi$, we have
$
1-(\frac {|\bar y| }\epsilon)^2\ge  C^{-1}$,
and $ 
C^{-1} \sup_{B_\epsilon}\hat u
\le \beta\epsilon^2 \le  C \inf_{ B_\epsilon}\hat u,
$
where $C$ is some positive constant depending
only on $(M,g)$.
Consequently,
$
A^{ \hat u}(\bar y)\le A^{\xi}(\bar y)\le C\beta^{ -\frac 4{n-2} }I.
$
This, together with the fact 
that $\lambda(A^{ \hat u}(\bar y))\in \Gamma\subset \Gamma_1$,
implies that
$
|\lambda (A^{ \hat u}(\bar y))|\le C \beta^{ -\frac 4{n-2} }.
$
Since $f(\lambda(A^{ \hat u}(\bar y)))=1$, we have, by (\ref{9new}),
that
$
\beta\le C\delta_1^{\frac {2-n}4},
$
where $C$ depends only on $(M,g)$.
Again by  (\ref{harnack}), we have
$$
\max_M u=\tilde u(\tilde x)\le C\hat u(0)
\le C\hat u(\bar y)=C\xi(\bar y)
\le C\beta\le C\delta_1^{\frac {2-n}4}.
$$
Namely, we have proved, for some positive constant
$C$ depending only on $(M,g)$, that
$
u\le C\delta_1^{\frac {2-n}4}$ on $M$.
Let $\bar x\in M$ be a maximum point of $u$, it was shown in
\cite{LL} that
$
f(u(\bar x)^{ -\frac 4{n-2} }\lambda(A_g(\bar x)))\le 1.
$
This, together with (\ref{delta2}), implies
$
\max_{M}u=u(\bar x)\ge \delta_2^{ \frac{n-2}4}.
$
Using  the upper bound of  $|\nabla_g \log u|$
on $M$, we have, for some positive constant $C$ depending only on
$(M,g)$, that
$
u\ge 
\frac 1C \max_M u\ge \frac 1C \delta_2^{ \frac{n-2}4}$, on $ M$.
Step 1 is established.

\noindent{\bf Step 2.}\  For some  positive constant  $C$
 depending only on $(M,g)$, $b$,  $\delta_1$ and $\delta_2$, we have
$
|\nabla_g^2u|\le C$ on $ M$.

$C^2$ estimates for
$(f,\Gamma)=(\sigma_k^{\frac 1k}, \Gamma_k)$
were obtained by Viaclovsky \cite{V1}.  The arguments
can be adapted in our situation.  Indeed, this is
equivalent to setting $\rho\equiv 1$ in the definition
of $G(x)$ in the proof of theorem 1.6 in
\cite{LL}, so that $G(x)$ is defined on
$M$, and Step 2 follows from the computation
there (with $h\equiv 1$) together with our $C^0$ and $C^1$ estimates
of $u$ and $u^{-1}$ obtained in Step 1.
Since $f$ is concave in $\Gamma$, and since
we have established $C^0$, $C^1$ and $C^2$ estimates of
$u$ and $u^{-1}$, higher derivative estimates of
$u$ and $u^{-1}$ in (\ref{aa7}) follow from  the interior
estimates of Evans 
and Krylov  together with the Schauder estimates.
Estimate (\ref{aa7}) has been established.

For the existence part of Theorem \ref{t1},
we only need to treat the case that  $(M,g)$ is not conformally diffeomorphic
to a standard sphere since it is obvious otherwise.
The following homotopy 
was introduced in \cite{LL}:
 For
$0\le t\le 1$,  let
$f_t(\lambda)=f\left(t\lambda+(1-t)\sigma_1(\lambda)e\right),
$ be defined on
$ \Gamma_t:=\{ \lambda\in \Bbb R^n\ |\
t\lambda+(1-t)\sigma_1(\lambda)e\in \Gamma\}.
$
We consider, for $0\le t\le 1$,
and for $\hat g=u^{ \frac 4{n-2}}g$,
\begin{equation}
f_t(\lambda(A_{\hat g}))=1,\quad\lambda(A_{\hat g})\in \Gamma_t,
\qquad \mbox{on}\ M.
\label{equationt}
\end{equation}

For $0\le t\le 1$, $(f_t, \Gamma_t)$ satisfies
(\ref{gamma0}), (\ref{gamma1}), (\ref{gamma1new}),
 (\ref{hypo0}), (\ref{aa5})
and (\ref{aa6}).  Moreover estimate (\ref{aa7}) holds
for solutions of  (\ref{equationt}), uniform in $0\le t\le 1$.
With this uniform estimates the degree argument in
\cite{LL} yields 
 a solution
$u$ of  (\ref{ghat2})  in $C^{4,\alpha}$.
By 
standard elliptic theories, $u\in C^\infty(M)$. 
 Theorem \ref{t1} is established.
 
Next we present our proof of
 Theorem~\ref{t2}.
By scaling, it is easy to see that we only need to prove the theorem for
$R=\delta=1$, which we assume below.
Let $u(\bar{x})=\max\limits_{\bar{B}_1}u$. 
As in the proof of theorem 1.8 in \cite{LL},we can find $\tilde{x}\in
B_{\frac 14}(\bar{x})$ such that
\begin{equation}
u(\tilde{x})\ge 2^{\frac{2-n}{2}}\sup_{B_{\sigma}(\tilde{x})}u,
\quad \mbox{and}\ \  \gamma:=u(\tilde{x})^{\frac{2}{n-2}}\sigma\ge\frac1 2
u(\bar{x})^{\frac{2}{n-2}},
\label{another2}
\end{equation}
where $\sigma=\frac12(1-|\tilde{x}-\bar{x}|)\le\frac12$.\newline

If 
$
\gamma\le 2^{n+8}n^4,
$
then
$
(\sup_{B_1}u)(\inf_{B_2}u)
\le u(\bar x)^2\le (2\gamma)^{ \frac {n-2}2 }
\le C(n),
$
and we are done.  So we always assume that
$
\gamma>  2^{n+8}n^4.
$
Let $\Gamma:=u(\tilde{x})^{\frac{2}{n-2}}\ge 2\gamma$, and consider 
$
w(y):=\frac{1}{u(\tilde{x})}u\Big(\tilde{x}+\frac{y}{u(\tilde{x})
^{\frac{2}{n-2}}}\Big),\ |y|<\Gamma.
$
By supharmonicity of $u$, 
\begin{equation}
\min_{\partial B_\Gamma}
w=\inf_{ B_\Gamma}w\ge \frac{1}{u(\tilde{x})}\min_{\partial B_2}u,
\quad 1=w(0)\ge 2^{\frac{2-n}{2}} \sup_{ B_\gamma}w.
\label{another3}
\end{equation}
By the conformal invariance of the equation satisfied by $u$,
$
F(A^w)=1$ on $B_\Gamma$.
Fix
$
r= 2^{n+6} n^4<\frac 1{4} \gamma.
$
For  $|x|<r$, 
 consider 
$$
w_{x,\lambda}(y):= (\frac \lambda{ |y-x| })^{n-2}
w(x+ \frac {\lambda^2(y-x) }{ |y-x|^2 }), \qquad
y\in B_\Gamma.
$$
By the conformal invariance of the equation, we have
$
F(A^{ w_{x,\lambda}})=1$ on
$B_\Gamma \setminus B_\lambda(x)$.
As in \cite{LL},  there exists $0<\lambda_x<r$ such that
we have
$$
w_{x,\lambda}(y)\le w(y),\qquad
\forall \ 0<\lambda<\lambda_x, \ y\in B_\Gamma\setminus B_\lambda(x),
$$
and
$$
w_{x,\lambda}(y)< w(y),\qquad
\forall \ 0<\lambda<\lambda_x, \ y\in \partial  B_\Gamma.
$$
By the moving sphere argument in \cite{LL},
we only need to consider the following two cases:

\noindent {\bf Case 1.}\  For some $|x|<r$ and some
 $\lambda\in (0,r)$,   $w_{x,\lambda}$ touches
$w$ on $\partial B_\Gamma$.

\noindent {\bf Case 2.}\  For all $|x|<r$ and all 
$\lambda\in (0,r)$, we have 
$$
w_{x,\lambda}(y)\le w(y),\quad \forall\  
|y-x|\ge\lambda, \ y\in B_\Gamma.
$$

In  Case 1, let $\lambda\in (0,r)$ be the smallest
number for which  $w_{x,\lambda}$ touches $w$
on $\partial B_\Gamma$.
By (\ref{another3}), we have, for some $|y_0|=\Gamma$,
$
u(\tilde{x})^{-1}\min_{\partial
B_2}u\le\min_{\partial B_\Gamma}w= w_{x,\lambda}(y_0).
$
Using (\ref{another3}), 
\[
w_{x,\lambda}(y_0)\le \Big(\frac{\lambda}{\Gamma-|x|}\Big)^{n-2}
\sup_{B_\gamma}w\le 
2^{\frac{n-2}{2}}\Big(\frac{\lambda}{\Gamma-|x|}\Big)^{n-2}\le
2^{\frac{n-2}{2}}\Big(\frac{r}{\Gamma-r}\Big)^{n-2}.
\]
Therefore
\[
\sigma^{\frac{n-2}{2}}u(\tilde{x})\min_{\partial B_2}u\le
2^{\frac{n-2}{2}}\sigma^{\frac{n-2}{2}}u(\tilde{x})^2
\Big(\frac{r}{\Gamma-r}\Big)^{n-2}.
\]
Since $4r<\gamma\le\frac{\Gamma}{2}$ and $\sigma\le\frac12$,
\begin{equation}
\sigma^{\frac{n-2}{2}}u(\tilde{x})\min_{\partial B_2}u\le
2^{\frac{n-2}{2}}\sigma^{\frac{n-2}{2}}u(\tilde{x})^2\frac{r^{n-2}}
{(\frac12\Gamma)^{n-2}}=
2^{\frac32(n-2)}\sigma^{\frac{n-2}{2}}r^{n-2}
\le C(n).
\label{another5}
\end{equation}
We deduce from  (\ref{another2}) and (\ref{another5})  that
$
(\sup_{B_1}u)(\inf_{B_2}u)\le 8^{n-2}r^{n-2}
\le C(n).
$

In Case 2, we have, by Lemma~\ref{l2}
and (\ref{another3}), that  
$$
|\nabla w(y)|\le 2(n-2)r^{-1}w(y)\le
(n-2)2^{ \frac {n}2 }r^{-1},\quad\forall |y|\le r.
$$
Let $\epsilon$ be the number such that
$
\xi(y):=\frac{1-\epsilon}{r}(r-|y|^2)$
satisfies
$
w\ge \xi$ on
$B_{\sqrt{r}}$ and
for some $|\bar y|<\sqrt{r}$,
$
w(\bar y)=\xi(\bar y)$.
Since $1=w(0)\ge \xi(0)=1-\epsilon$ and 
$\xi(\bar y)>0$, we have $0\le \epsilon<1$.

By the estimates of $|\nabla w|$ and the mean value theorem,
$
|w(y)-1|\le (n-2)2^{ \frac {n}2} r^{-\frac 12},
$ for all $|y|\le \sqrt{r}$.
So 
$
\frac 12\le
1-  (n-2)2^{ \frac {n}2} r^{-\frac 12}\le w(\bar y)=\xi(\bar y)\le 
1-\epsilon,
$
and therefore
$
0\le \epsilon\le  (n-2)2^{ \frac {n}2} r^{-\frac 12}.
$
Clearly,
$$
\nabla w(\bar y)=\nabla \xi(\bar y),
\ \ |\nabla \xi(\bar y)|\le \frac 2{\sqrt{r}},
\ \ D^2w(\bar y)\ge D^2\xi(\bar y)= -2(1-\epsilon)r^{-1}I.
$$
It follows that
\[
A^w(\bar y)\le
A^{\xi}(\bar y)\le
\frac{(10n+4)}{ (n-2)^2 } 2^{ \frac {2n}{n-2} } r^{-1}I.
\]
Since $F(A^w(\bar y))=1$, we have, by (\ref{F3}) (recall
that $\delta=1$),
$ \frac{(10n+4)}{ (n-2)^2 } 2^{ \frac {2n}{n-2} } r^{-1}
\ge 1$, violating the choice of $r$.
Thus we have shown that Case 2 can never occur.
Theorem~\ref{t2} is established.

The results in this note
have been presented by the second
author at his 45-minute 
invited talk at  ICM 2002  in August 2002 in Beijing.
The results have also been presented by the second author
in a colloquium talk at Northwestern University on
September 27, 2002, in the Geometric Analysis
seminar at Princeton University 
on October 18, 2002, in a mini-course in late October 2002
 at Universit\`a di Milano.
On December 2 2002, the second author was informed 
by P. Guan that he, in collaboration with C.S. Lin and G. Wang,
has obtained some  related results.

\end{document}